\documentstyle[12pt]{article}
\makeatletter
\begin{document}
\thispagestyle{empty}
\pagestyle{empty}
\parskip = 10pt
\newcommand{\namelistlabel}[1]{\mbox{#1}\hfil}
\newenvironment{namelist}[1]{%
\begin{list}{}
{
\let\makelabel\namelistlabel
\settowidth{\labelwidth}{#1}
\setlength{\leftmargin}{1.1\labelwidth}
}
}{%
\end{list}}
\def\theequation{\thesection.\arabic{equation}}
\newtheorem{theorem}{\bf Theorem}[thetheorem]
\newtheorem{corollory}{\bf Corollory}[thecorollory]
\newtheorem{remark}{\bf Remark}[theremark]
\newtheorem{lemma}{\bf Lemma}[thelemma]
\def\thetheorem{\arabic{theorem}}
\def\theremark{\thesection.\arabic{remark}}
\def\thecorollory{\thesection.\arabic{corollory}}
\def\thelemma{\thesection.\arabic{lemma}}
\newcommand{\nc}{\newcommand}
\newcommand{\bsp}{\begin{sloppypar}}
\newcommand{\esp}{\end{sloppypar}}
\newcommand{\be}{\begin{equation}}
\newcommand{\ee}{\end{equation}}
\newcommand{\beanno}{\begin{eqnarray*}}
\newcommand{\inp}[2]{\left( {#1} ,\,{#2} \right)}
\newcommand{\dip}[2]{\left< {#1} ,\,{#2} \right>}
\newcommand{\disn}[1]{\|{#1}\|_h}
\newcommand{\pax}[1]{\frac{\partial{#1}}{\partial x}}
\newcommand{\tpar}[1]{\frac{\partial{#1}}{\partial t}}
\newcommand{\xpax}[2]{\frac{\partial^{#1}{#2}}{\partial x^{#1}}}
\newcommand{\pat}[2]{\frac{\partial^{#1}{#2}}{\partial t^{#1}}}
\newcommand{\ntpa}[2]{{\|\frac{\partial{#1}}{\partial t}\|}_{#2}}
\newcommand{\xpat}[2]{\frac{\partial^{#1}{#2}}{\partial t \partial x}}
\newcommand{\npat}[3]{{\|\frac{\partial^{#1}{#2}}{\partial t^{#1}}\|}_{#3}}
\newcommand{\eeanno}{\end{eqnarray*}}
\newcommand{\bea}{\begin{eqnarray}}
\newcommand{\eea}{\end{eqnarray}}
\newcommand{\ba}{\begin{array}}
\newcommand{\ea}{\end{array}}
\newcommand{\nno}{\nonumber}
\newcommand{\dou}{\partial}
\newcommand{\bc}{\begin{center}}
\newcommand{\ec}{\end{center}}
\newcommand{\bb}{\mbox{\hspace{.25cm}}}
\nc{\benu}{\begin{enumerate}}
\nc{\eenu}{\end{enumerate}}
\nc{\bth}{\begin{theorem}}
\nc{\eth}{\end{theorem}}
\nc{\bpr}{\begin{prop}}
\nc{\epr}{\end{prop}}
\nc{\blem}{\begin{lemma}}
\nc{\elem}{\end{lemma}}
\nc{\bcor}{\begin{corollary}}
\nc{\ecor}{\end{corollary}}
\newcommand{\R}{I\!\!\!R}
\nc{\Hy}{I\!\!\!H}
\nc{\C}{I\!\!\!C}
\newcommand{\la}{\lambda}
\newcommand{\s}{\sinh}
\newcommand{\co}{\cosh}
\newcommand{\vl}{V_{\la}}
\nc{\ga}{\gamma}
\nc{\vg}{V_{\ga}}
\nc{\T}{\Theta}
\nc{\f}{\frac}
\nc{\rw}{\rightarrow}
\nc{\om}{\omega}
\nc{\Om}{\Omega}
\nc{\al}{\alpha}
\nc{\qed}{\hfill \rule{2.5mm}{2.5mm}}
\setcounter{section}{0}
\title{Finite dimensional imbeddings of harmonic spaces}
\author{K. Ramachandran\thanks{supported by National Board for 
Higher Mathematics, DAE, INDIA} \,and A. Ranjan}
\maketitle
\begin{abstract}
In a noncompact harmonic manifold $M$ we establish finite dimensionality of
the eigenspaces $V_{\lambda}$ generated by 
radial eigenfunctions of the form $\cosh r + c$. As a consequence, for such 
harmonic manifolds, we give an isometric imbedding of $M$ into
$(V_{\lambda},B)$, where $B$ is a nondegenerate symmetric bilinear
indefinite form on $V_{\lambda}$ (analogous to the imbedding of
the real hyperbolic space
$I\!\!\!H^{n}$ into $\R^{n+1}$ with the indefinite form
$Q(x,x) = -x_0^2 + \sum x_i^2$). Finally we give certain conditions under
which $M$ is symmetric.
\end{abstract}
\section{Introduction}
Let $(M,g)$ be a riemannian manifold. $(M,g)$ is said to be harmonic  
if it's density function $\omega_{p}=\sqrt{| \mbox{det}(g_{ij})|}$
in a normal coordinate neighbourhood around each point $p$ depends 
only on the geodesic distance $r(p,.)$. Let $(r_{p},\phi)$ be polar
coordinates around $p$, where $r_{p}(q) = r(p,q)$ is the geodesic 
distance between $p$ and $q$ and $\phi$ represents a point on the 
unit sphere in the tangent space $T_p(M)$ at $p$. In this system
the density is 
$$\Theta_p\,=\,r^{n-1}\,\omega_{p}$$
A riemannian manifold is said to be globally harmonic if the density
$\Theta_p$ is a radial or a spherically symmetric function around each
point $p \in M$ i.e it depends only on the vaiable $r$ and can thus be
written as $\Theta_p(r)$. The function $\T_p(r)$ does not depend on
the point $p$, since $\om_p(q) = \om_q(p)$ i.e the function $\om_p$
itself is independent of the point $p$.

Let $(M,g)$ be a harmonic manifold and let $\Delta$ be the laplacian on
$M$. Consider a nonzero eigenvalue $\lambda$ and the corresponding eigen
space $V_{\lambda}$. If $M$ is compact, $V_{\lambda}$ is known to be
finite dimensional. In this case Besse \cite{Bes} constructed isometric
imbedding of   
$M$ into $V_{\lambda}$. This immersion was crucially used by Szabo 
\cite{Sza} in his proof of the Lichnerowicz conjecture, which asserts
that {\it Harmonic manifolds are symmetric}. Szabo proved the conjecture 
for compact harmonic manifolds with finite fundamental group. In his paper
he generalised Besse's imbedding theorem to give isometric immersions of
any harmonic manifold $M$ into $L^2(M,g)$, an infinite dimensional space.
In particular he recovers Besse's imbedding for compact harmonic manifold.  
The eigen
spaces $V_{\lambda}$ of noncompact manifolds, in particular noncompact
harmonic manifolds need not be finite dimensional. Nevertheless the
problem of imbedding $M$ into a finite dimensional vector space albeit
with an indefinite metric is important and has not been addressed
(analogous to the imbedding of the real
hyperbolic space $I\!\!\!H^{n}$ into $\R^{n+1}$ with the
indefinite metric $Q(x,x) = -x_0^2 + \sum x_i^2$ ). We prove a partial 
result in this direction.
\begin{theorem}
Let $(M,g)$ be non-compact harmonic manifold. Let $0 \neq \lambda$ be an
eigenvalue of the Laplacian $\Delta$. Assume that all the averaged
eigenfunctions of $\Delta$ corresponding to $\lambda$
are of the form $\cosh r + c$, where c is a constant and
$r$ is the distance function from the point of averaging. Let
$V_{\lambda}$ be the eigen subspace generated by radial eigen functions.
Then
\begin{enumerate}
\item  $V_{\lambda}$ is finite dimensional and
\item there exists a symmetric nondegenerate bilinear form $B$ on
$V_{\lambda}$ and an isometric 
imbedding of $M$ into $( V_{\lambda}\,,\,-B )$.
\end{enumerate}
Moreover in this case the volume density is\\
$$\Theta(r)\,=\,2^{n-1}\,(\sinh \frac{r}{2})^{n-1}\,(\cosh
\frac{r}{2})^{n-1-\frac{4}{3}(k+n-1)}.$$ 
Also $c$ and $\lambda$ are given by
$$c\,=\,\frac{k+n-1}{\frac{n}{2}+1+k}$$
and $$\lambda\,=\,-\f{2}{3}(\frac{n}{2}+1+k)$$
where $Ricci(M) = -k, k > 0$. 
Finally $f_p = \co r + c$ has no zeroes on $M$.
\end{theorem}
Damek and Ricci \cite{DR} have given a family of noncompact harmonic 
spaces which are not symmetric, thus disproving the Lichnerowicz
conjecture for the noncompact case. These spaces are called the NA spaces.
Though these spaces are non-symmetric, they have an eigen function of the
form $\cosh r + C$. These are the only known examples of non-symmetric
harmonic spaces. Thus our theorem covers all the known cases of harmonic 
manifolds. Nevertheless we strongly believe that our theorem is true for
all harmonic spaces.
\section{The imbedding of $M$ into $(V_{\lambda},-B)$}
Before proving the theorem we make a simple but useful observation.\\
{\bf Observation :} Let $V$ be a vector space, $S \subseteq V$ a set 
such that $V = span S$. Let 
$\tilde{B} : S \times V \rightarrow  I\!\!\!R$
be a map such that\\
1. $\tilde{B}$ is linear in the second variable, i.e, $\tilde{B}
   (s , .)$ is linear on $V$ for all $s$ in $S$.\\
2. $\tilde{B} (s , t) = \tilde{B} (t , s)\,\,\forall\, s , t \in S.$\\
then there exists a unique bilinear form $B : V \times V \rightarrow 
I\!\!\!R $ which extends $\tilde{B}$.\\
{\bf Proof :} Define 
$$
B\,(\sum_i\,a_i\,s_i\,,\,v)\,=\,\sum_i\,a_i\,\tilde{B}\,(s_i\,,\,v)$$
\qed

\noindent
{\bf Proof of the Theorem :}\\
Let $x$ be any point of $M$ and let $B_x$ be the unit ball around $x$. 
Let $C^{\infty} (B_x)$ be the space of smooth functions on $B_x$ with 
the {\it sup} norm. 

Consider the restriction map 
$$\Psi\,:\,V_{\lambda}\,\,\rightarrow\,\,C^{\infty}\,(B_x)\,;
\,f\,\mapsto\,f|_{B_x}$$
Due to analyticity $\Psi$ is injective. We will show that $Im \Psi$ is finite
dimensional. This will prove that $V_{\lambda}$ is finite dimensional.

Let $f$ be an eigenfunction for the eigen value $\lambda$ and let $f_p$
denote the averaged eigen function around $p$. We denote by
$C\,=\,V_{\lambda}^+$ the cone generated by $f_p , p \in M$, i.e
$V_{\lambda}^+ = \{ \sum_i a_p f_p : a_p \geq 0 , p \in M \}$ and only 
finite sums are taken. 

\noindent
{\bf claim :}  $V_{\lambda} \cap B_x$ is finite dimensional, i.e $Im \Psi$
is finite dimensional.\\
{\bf proof :} Let $ h = \sum_{1}^{k} a_{p_i} f_{p_i} \in C$. 
Since $f_{p_i} = \cosh d(p_i,.)$  
$$ \nabla\,h\,=\,\sum_{1}^{k}\,a_{p_i}\,\sinh\,d(p_i,.)$$
Hence 
\begin{eqnarray*}
|\,\nabla\,h|^2\,&=& (\,\sum\,a_{p_i}\,\sinh\,d(p_i,.)\,)^2 \\
                 &\leq&\,(\,\sum\,a_{p_i}\,cosh\,d(p_i,.))^2 \\
                 &=&\,|h|^2
\end{eqnarray*}
We consider 
the space $X = C^{\infty} (B_x) \cap V_{\lambda}^+$. Since 
$|\nabla h|_{B_x}| \leq |\nabla h|$, the above inequality
shows that for $h \in X$, $|\nabla h| \leq |h|$ holds. This shows the
following : If $W$ is a open ball in $V_{\lambda}^+$ such that 
$W \cap C^{\infty} (B_x)$ is a bounded family then it is a pointwise
bounded equicontinuous family in $X$. Hence by the Arzela-Ascoli theorem
it is relatively compact in $X$.

Now we argue as follows: let $U$ be the open ball around the origin
in $V_{\lambda}$. Choose a $g$ in  $ V_{\lambda}$ such
that $U\,\subset\, C\,-\,g$ i.e $U-g\,\subset\,C$. The above argument
shows that $U-g\,\cap\,C^{\infty}(B_x)$ is relatively compact in
$V_{\lambda} \cap C^{\infty}(B_x)$, hence $V_{\lambda} \cap
C^{\infty}(B_x)$ is finte dimensional and the claim is proved.

Let $S = \{ f_p : p \in M \}$ then $V_{\lambda} = Span S$. Define 
$$\tilde{B}\,\,:\,\,S\,\times\,V_{\lambda}\,\rightarrow\,\R\,\,{\rm by}$$
$$\tilde{B}\,(\,f_p\,,\,\sum\,\alpha_q\,f_q\,)\,=\,
\sum\,\alpha_q\,f_q\,(p)$$
Then by the above observation $\exists !\,
B\,:\,V_{\lambda}\,\times\,V_{\lambda}\,\rightarrow\,\R$, which extends
$\tilde{B}$. \\
{\bf claim:} $B$ is symmetric and nondegenerate.\\
{\bf proof:} Clearly $B$ is symmetric. Nondegeneracy of $B$ can be 
established as follows. Let $h \in KerB$, then 
$$ B\,(h,f_q)\,=\,0\,\forall\,q \in M\,i.e,\, h(q)\,=\,0\,\forall\,q \in M$$  
Now define the map\\
$$ \Phi\,\,:\,\,M\,\,\rightarrow\,\,(\,V_{\lambda}\,,\,-B\,)\,\,{\rm by}$$
$$ \Phi\,(p)\,:=\,\Phi_p\,=\,\cosh\,d\,(\,p\,,\,.)$$
Note that 
$\Phi(M)$ is contained in $\{ B(h,h) = 1\}$ for let
$\gamma(t)$ be a geodesic in $M$ then 
$$B\,(\,\Phi_{\gamma(t)}\,,\,\Phi_{\gamma(t)}\,)\,=\,\Phi_{\gamma(t)}\,
\gamma(t)\,=\,1$$
Further 
$$B\,(\,\Phi_{\gamma(t)}\,,\,\Phi_{\gamma(t)}\,)\,=\,1$$
gives on twice differentiating,
$$B\,(\,\Phi'_{\gamma(t)}\,,\,\Phi'_{\gamma(t)}\,)\,|_{t\,=\,0}\,+\,
B\,(\,\Phi''_{\gamma(t)}\,,\,\Phi_{\gamma(t)}\,)\,|_{t\,=\,0}\,=\,0$$
$${\rm
i.e}\,\,B\,(\,\Phi'_{\gamma(t)}\,,\,\Phi'_{\gamma(t)}\,)\,|_{t\,=\,0}\,+\,
B\,(\,\Phi''_{\gamma(t)}\,|_{t\,=\,0},\,\Phi_{\gamma(0)}\,)\,=\,0$$
But, since
$\Phi_{\gamma(t)} = \cosh d\,(\gamma(t))$, 
$\Phi''_{\gamma(t)} = \cosh d (\gamma(t),.)$ and hence
$$B\,(\,\Phi''_{\gamma(t)}\,|_{t\,=\,0},\,\Phi_{\gamma(0)}\,)\,=\,1 $$
Therefore 
$$B\,(\,\Phi'_{\gamma(t)}\,,\,\Phi'_{\gamma(t)}\,)\,|_{t\,=\,0}\,=\,-1$$
This shows that $\Phi$ gives an isometric immersion of $M$ into 
$( V_{\lambda} , -B )$. Finally, since $\cosh r$ is a monotone function,
$\Phi_p \neq \Phi_q$ for $p \neq q$. Thus we get an imbedding of $M$
into $(\vl , -B)$. This imbedding is nondegenerate since the image of $M$ 
is not contained in a lower dimensional subspace of $\vl$.

We now compute the densities of these harmonic manifolds. 
Let $(r,\phi) $ be polar coordinates around any point $p$ of $M$ and 
$f$ be any smooth function on $M$. 
The Laplacian $\Delta$ in these coordinates is given by
$$\Delta f\,=\,\Delta_{r}f\,-f''\,-\,\frac{\Theta'}{\Theta}f'$$
Where $\Delta_{r}$ is the laplacian on the distance sphere $S_r$ and 
the primes denote derivatives with respect to $r$. Therefore the action of
$\Delta$ on a radial function $f$ (around $p$) is
$$\Delta f\,=\,-f''\,-\,\frac{\Theta'}{\Theta}f'$$
Hence if $f$ is a radial eigen function for the eigen value $\lambda$, 
then $f$ satisfies 
$$f''\,+\,\frac{\Theta'}{\Theta}f'\,+\,\lambda f\,=\,0$$
By our assumption $f(r) = \cosh r + c$. Hence we get
$$\frac{\Theta'}{\Theta}(r)\sinh r\,=\,-(1+\lambda)\cosh r\,-\,c\lambda$$
which gives 
$$\T(r)\,=\,2^{-(1+\la)}\,(\s(\f{r}{2}))^{-(1+\la+c\la)}\,(\co
(\f{r}{2}))^{-(1+\la-c\la)} $$ 
But as $r \rw 0,\,\T(r) \rw (\s (\f{r}{2}))^{n-1}$, hence we get
$$1+\la+c\la\,=\,1-n,\, {\rm and}\,\,\T(r)\,=\,2^{(n-1)}\,
(\s (\f{r}{2}))^{n-1}\,(\co (\f{r}{2}))^{(n-1+2c\la)}$$  
We now use the ledgers formula \cite{Bes} (see pp. 161), which gives
$$-\f{1}{3} Ricci\,=\,\om''(r)|_{r=0}$$
This shows that harmonic manifolds are Einstein, i.e they have 
constant ricci curvature. Let $Ricci(M) = -k, k>0$.
Now $$\T(r)\,=\,r^{n-1}\om(r)$$
gives $$\f{\T'(r)}{\T(r)}\,=\,\f{n-1}{r}\,+\,\f{\om'(r)}{\om(r)}$$
Using the expression for $\T$ obtained above, we get
$$\f{\om'(r)}{\om(r)}\,=\,-(1+\la)\coth r\,-\,\f{c\la}{\s
r}\,-\,\f{n-1}{r}$$ 
Differentiating twice and after simplification, we get
$$-\f{1}{3} Ricci\,=\,\f{n-1}{3}\,+\,\f{c\la}{2}$$
which gives
$$c\la\,=\,-\f{2}{3}(k+n-1), {\rm but}\,1+\la+c\la\,=\,1-n $$
hence $$ \la\,=\,-\f{2}{3}(1+k+\f{n}{2})\,{\rm and}\,
         c\,=\,\f{k+n-1}{1+k+\f{n}{2}}$$
Finally the above computations show that
\beanno 
f\,&=&\,\co r\,+\,c\\
   &=&\,\co r\,-\,1\,+\,\f{3n/2}{1+k+\f{n}{2}}\\
   &=&\,2\s^2 (\f{r}{2})\,+\,\f{3n/2}{1+k+\f{n}{2}}\\
   &\geq&\,\f{3n/2}{1+k+\f{n}{2}}
\eeanno
Thus $f$ has no zeroes on $M$. This completes the proof of  the theorem.
\qed
 
\section{Symmetry of $M$}
In ths section we prove the symmetry of $M$ under some assumptions.
Two cases arise.\\
{\bf Case 1. c = 0}\\
In this case $M$ is symmetric. The sturm liouville equation 
$$f''\,+\,\frac{\T'}{\T}f\,+\,\la f\,=\,0$$
shows that 
$$\T(r)\,=\,\s^{n-1}r$$
A result from \cite{RR} shows that $M$ is isometric to the real hyperbolic
space with constant curvature -1.\\
{\bf Case 2. c $\neq$ 0}\\
Let $\vg$ be the subspace of $\vl$ generated by $\{f_{\ga(t)} : t \in \R\}$.
\begin{lemma}
Let $M$ satisfy the hypothesis of the above Theorem. Assume in addition
that $B$ restricted to the subspace $\vg$ is nondegenerate. Then $M$ is
symmetric.
\end{lemma}
{\bf Proof.} Fix an imbedding $\Phi : M \rw \vl ; \Phi(p) = f_p$. Let 
$\ga(t)$ be a geodesic on $M$. Then\\
$$B(f_{\ga(t)},f_{\ga(0)})\,=\,\co t\,+\,c$$
Differentiating we get
$$B(f'''_{\ga(t)}-f'_{\ga(t)},f_{\ga(0)})\,=\,0$$
Since $\ga(0)$ is an arbitrary point of$M$, $f'''_{\ga(t)}-f'_{\ga(t)} \in
\vg$ and by our assumption $B$ is nondegenerate on $\vg$, the geodesic
$\ga$ satisfies 
$$f'''_{\ga(t)}-f'_{\ga(t)}\,=\,0$$
Take $p \in M$ and let $N_p$ be the orthogonal complement, w.r.t $B$ 
of the tangent space $T_p(M)$ to $M$. Define
$$\tau_p\,:\,\vl\,\rw\,\vl$$
be the reflection w.r.t the subspace $N_p$. $\tau_p$ is an isometry 
of $\vl$, it fixes $p$ and reverses the geodesics through $p$.
Moreover it also leaves $M$ invariant. hence it induces an isometry 
on $M$ which is obviously the geodesic involution of $M$. Thus $M$ 
is a symmetric space.
\qed
\newpage
The next lemma shows that if $\vg$ is 3-dimensional then $B$ is
nondegenerate on it.
\begin{lemma}
Let $M$ satisfy the hypothesis of the theorem. Let $\ga$ be a geodesic
such that $\vg$ is 3-dimensional, then $B$ is nondegenerate on $\vg$.
\end{lemma}
{\bf Proof.} Let $s_i = -1,0,1$ be such that $\{f_{\ga(s_i)}\}_1^3$ 
generate $\vg$. Let $g \in \vg$ satisfy $B(g,f_(\ga(t))) = 0\,\,\forall\,t$. 
Put 
\beanno
g\,&=&\,\sum_1^3\,a_i\,f_{\ga(s_i)}\\
   &=&\,\sum_1^3\,a_i\,\co d(\ga(s_i),.)\,+\,c\sum_1^3\,a_i
\eeanno
Now $B(g,f_{\ga(t)}) = 0$ gives $g(\ga(t) = 0 \,\forall t$, i.e
$$\sum\,a_i\,c\,+\,\sum\,a_i\,\co(s_i-t)\,=\,0\,\,\,\forall\,t$$
Dividing by $\co t$ and letting $t \rw \infty$ we see that 
$$\sum\,a_i\,=\,0$$
hence 
$$\sum\,a_i\,\co (s_i-t)\,=\,0\,\,\,\forall\,t$$
Expanding and equating the coefficients of $\co t$ and $\s t$ to 0
we get
$$\sum\,a_i\,\co s_i\,=\,0\,=\,\sum\,a_i\,\s s_i$$
Substituting the vaues of $s_i$ and using the equation $\sum a_i = 0$ 
we obtain a system of three equations 
\beanno
a_0\,\co 1\,+\,a_1\,+\,a_2\,\co 1\,&=&\,0\\
-a_0\,\s 1\,+\,0\,+\,a_2\,\s 1\,&=&\,0\\
a_0\,+\,a_1\,+\,a_2\,&=&\,0
\eeanno
But the coefficient matrix of the above system is nonsingular, hence 
$a_i = 0 \,\,\forall i$, i.e $g \equiv 0$ or $B$ is nondegenerate.
\qed\\
Combining the above two lemmas one sees that $M$ is symmetric if $\vg$
is 3-dimensional.
\newpage
\section{Summary}
For manifolds which satisfy the hypothesis of the theorem the density
$\T(r)$ is given by
$$\Theta(r)\,=\,2^{n-1}\,(\sinh \frac{r}{2})^{n-1}\,(\cosh
\frac{r}{2})^{n-1-\frac{4}{3}(k+n-1)}.$$ 
Put $b = n-1-\frac{4}{3}(k+n-1)$. By the Bishop-Gromov comparison
theorem $b$ satisfies $0 \leq b \leq n-1$.\\
\begin{enumerate}
\item If $b = n-1$, $\T(r) = \s^{n-1} r$ and hence $M$ is isometric to
the real hyperbolic space. In this case $c = 0$.  
\item If $b = 0$, $\T(r) = 2^{n-1} (\sinh \frac{r}{2})^{n-1}$. Again $M$
is isometric to the real hyperbolic space. But in this
case $c \neq 0$ and the eigen function $f_p = \co r + c$ becomes the
eigen function for the next eigen value.
\item All the other harmonic spaces which satisfy the hypothesis of the 
theorem, in particular the $NA$ spaces have $0 < b < n-1$.
\end{enumerate}
{\bf Question} Is $b$ an integer?, in that case are all harmonic spaces 
homogeneous?

\newpage
\bibliographystyle{plain}

\begin{thebibliography}{99}
\bibitem{Bes} A.L.Besse, {\em Manifolds all of whose geodesics are
closed}, Springer, Berlin, 1978. 
\bibitem{DR} E.Damek and F.Ricci, {\em A class of nonsymmetric harmonic
riemmanian spaces}, Bulletin AMS, 27(1992), pp. 139-142.
\bibitem{RR} K. Ramachandran and A. Ranjan, {\em Harmonic manifolds with
some specific volume densities}, preprint.
\bibitem{Sza} Z.I.Szabo, {\em The Lichnerowicz conjecture on Harmonic
Manifolds},  J. Differential Geometry, 31(1990), pp. 1-28. 
\bibitem{Sza1} Z.I.Szabo, {\em Spectral theory for operator families on
Riemannian Manifolds}, Proceedings of symposia in Pure Math, 54 (1993), 
part 3, pp 615-665.
\end{thebibliography}

\noindent
Department of Mathematics,\\
I.I.T Powai, Mumbai 400 076,\\
India.\\
email: kram@ganit.math.iitb.ernet.in\\
\hspace*{1.25cm}aranjan@ganit.math.iitb.ernet.in\\
\end{document}